\documentclass[10pt,conference]{ieeeconf}
\usepackage{amsmath,amssymb,latexsym,epsfig,psfrag,graphicx}
\newtheorem{prp}{Proposition}
\newtheorem{thm}[prp]{Theorem}

\newenvironment{pf}{\noindent{\it Proof. }}{\hfill$\Box$\smallbreak}
\newenvironment{pf*}[1]{\smallbreak\noindent{\it #1}}{\hfill$\Box$\smallbreak}
\newcounter{definition}

\newcounter{remark}
\newenvironment{rmk}{\addtocounter{remark}{1}\smallbreak\noindent
  {\em Remark \theremark.}}{\smallbreak}
\newcounter{example}

\newenvironment{ex*}[1]{\addtocounter{example}{1}\smallbreak\noindent
	{\bf Example \theexample{} --- {\bf #1}}}{\hfill$\Box$\smallbreak}
\newcommand{\realR}{\mathbb{R}}

\newcommand{\averageE}{\mathbb{E}}
\newcommand{\itemi}{($i$)}
\newcommand{\ii}{($ii$)}

\DeclareMathOperator{\conv}{conv}
\DeclareMathOperator{\diag}{diag}

\DeclareMathOperator*{\minnu}{\hbox{$\min_\nu$}}
\DeclareMathOperator*{\minv}{\hbox{$\min_v$}}

\DeclareMathOperator*{\maxxi}{\hbox{$\max_\xi$}}

\DeclareMathOperator{\vect}{vec}

\begin{document}
\title{Minimax Optimal Adaptive Control for Systems on Cones}
  \author{Anders Rantzer
  \thanks{The author is affiliated with Automatic Control LTH, Lund
    University, Box 118, SE-221 00 Lund, Sweden. He is a member of the Excellence Center ELLIIT and Wallenberg AI, Autonomous Systems and Software Program (WASP). Support was received from the European Research Council (Advanced Grant 834142) }}
\maketitle

\begin{abstract}
	The theory of optimal control on positive cones has recently
	identified several new problem classes where the Bellman equation can be solved explicitly, in analogy with classical linear quadratic control. In this
	paper, the idea is extended to minimax adaptive control, yielding exact solutions to instances of the Bellman equation for dual control. In particular, this allows for optimization of the fundamental tradeoff between exploration and
	exploitation.
\end{abstract}

\section{Introduction}
\label{sec:intro}
This paper is about adaptive control with worst-case models for disturbances and uncertain parameters. More precisely, the following problem is addressed: 

\smallskip

\textbf{General Problem:}
{\it Consider proper cones $\mathcal{U}\subseteq\realR^{N_y+N_v}$, $\mathcal{X}\subseteq\realR^{N_y}$, $\mathcal{V}\subseteq\mathcal{U}\times\mathcal{X}$ and $\mathcal{S}\subseteq\realR^{2N_y+N_v}$. Suppose that for every $y\in\mathcal{X}$ there exists $v$ with $(y,v)\in\mathcal{U}$. Find, if possible, a control law 
	$v_t:=\mu_t(y_0,\ldots,y_t,v_0,\ldots,v_{t-1})$ 
	such that $(y_t,v_t)\in\mathcal{U}$ and 
	\begin{align}
		\max_{s\in\mathcal{S}}\sum_{t=0}^{T-1}s^\top\begin{bmatrix}
			y_t\\v_t\\y_{t+1}
		\end{bmatrix}&\le0
		\label{eqn:generalcost}
	\end{align}
	for all $y_0,y_1,\ldots\in\mathcal{X}$ with $(y_t,u_t,y_{t+1})\in\mathcal{V}$ for all $t$ and $y_0=0$.}

\smallskip

The problem can be addressed as a dynamic game problem, with a minimizing player selecting $v$ and a maximizing player selecting $y$ and $s$. It will be solved by dynamic programming while exploiting the special structure. The vector $s$ encapsulates both system dynamics and cost function, so the set $\mathcal{S}$ represents a priori information about these. Moreover, it will be verified that a solution, if it exists, can always be written on the form
\begin{align*}
	v_t&=\eta\left(y_t,\sum_{t=0}^{t-1}(y_\tau,v_\tau,y_{\tau+1})\right),
\end{align*}
where the sum is a compression of all data available to the controller at time $t$. Whithout any data, all $s\in\mathcal{S}$ carry equal weight, but the presence of past data $\sum_{t=0}^{t-1}(y_\tau,v_\tau,y_{\tau+1})$ makes it possible to focus attention on a smaller set of $s$-vectors in reducing the cost. Hence the optimal controller will make a tradeoff between "exploration and exploitation".

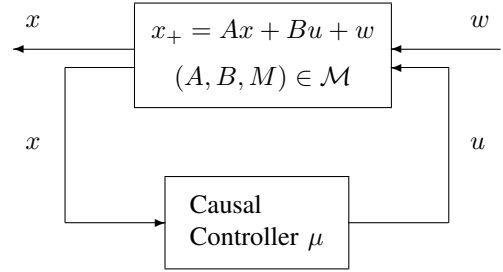
\begin{figure}
	\begin{center}
		\setlength{\unitlength}{.0081mm}%
		\begin{picture}(9366,4600)(1318,-8461)
			\put(4501,-8461){\framebox(3000,1500){\begin{tabular}{ll}Causal\\ Controller $\mu$
			\end{tabular}}}
			\put(4000,-5811){\framebox(4200,1700){\begin{tabular}{c}
						$x_+=Ax+Bu+w$\\[2mm]
						$(A,B,M)\in\mathcal{M}$
			\end{tabular}}}
			\put(4000,-5161){\line(-1, 0){1150}}
			\put(2851,-5161){\line( 0,-1){2550}}
			\put(2851,-7711){\vector( 1, 0){1650}}
			\put(4000,-4861){\vector(-1, 0){2000}}
			\put(10001,-4861){\vector(-1, 0){1800}}
			\put(7501,-7711){\line( 1, 0){1650}}
			\put(9151,-7711){\line( 0, 1){2550}}
			\put(9151,-5161){\vector(-1, 0){950}}
			\put(2200,-6500){$x$}%
			\put(9500,-6500){$u$}%
			\put(2200,-4500){$x$}%
			\put(9500,-4500){$w$}%
		\end{picture}%
	\end{center}
	\caption{We want a feedback controller that works for all system parameters within the given bounds.  Nonlinear adaptive controllers can do much better by estimating $(A,B)$ and use the estimate for control. The objective of this paper is to derive such controllers using a dynamic game formulation.
	}
	\label{fig:blockdiag}
\end{figure}

To clarify the connection to adaptive control, it is instructive to consider linear systems
\begin{align}
	x_{t+1}&=Ax_t+Bu_t+w_t&x_0&=0
\label{eqn:state}
\end{align}
with the objective that
\begin{align}
	\averageE\sum_{t=0}^{T-1}{\footnotesize\begin{bmatrix}x_t\\u_t\end{bmatrix}^\top
		M\begin{bmatrix}x_t\\u_t\end{bmatrix}}
	\le \gamma^2\sum_{t=0}^{T-1}|w_t|^2
\label{eqn:gainbound}
\end{align}
for all $w_0,w_1,\ldots$ and all $(A,B,M)$ in a given set $\mathcal{M}$. See Figure~\ref{fig:blockdiag}.
Eliminating $w_t$ from \eqref{eqn:gainbound} using \eqref{eqn:state} this problem can be turned into the following form:
\medskip

\textbf{Special problem --- Linear Quadratic Control:}
	{\it Given $\gamma>0$ and a set of matrix triples $(A,B,M)$, find a control law $u_t=\mu_t(x_0,\ldots,x_t,w_t)$ such that
\begin{align*}
	\averageE\sum_{t=0}^{T-1}\left({\footnotesize\begin{bmatrix}x_t\\u_t\end{bmatrix}^\top
		M\begin{bmatrix}x_t\\u_t\end{bmatrix}}-\gamma^2|Ax_t+Bu_t-x_{t+1}|^2\right)\le0
\end{align*}
	for all sequences $x_1,x_2,\ldots\in\realR^n$ and all $(A,B,M)$.}

\section{Notation}
The set of $n\times m$ matrices with real coefficients is denoted $\realR^{n\times m}$. The transpose of a matrix $A$ is denoted $A^\top$. For a symmetric matrix $A\in\realR^{n\times n}$, we write $A\succ0$ to say that $A$ is positive definite, while $A\succeq0$ means positive semi-definite. 
For $A, B\in\realR^{n\times m}$, the expression $\langle A,B\rangle$ denotes the trace of $A^\top B$. 
Given $x\in\realR^n$ and $A\in\realR^{n\times n}$, the notation $|x|^2_A$ means $x^\top Ax$. Similarly, given $B\in\realR^{n\times m}$ and $A\in\realR^{n\times n}$, the trace of $B^\top AB$ is denoted $\|B\|^2_A$. For two vectors $x\in\realR^n$, $u\in\realR^m$, the concatenation {\footnotesize$\begin{bmatrix}x\\u\end{bmatrix}$} will generally be denoted $(x,u)$. In particular, $|(x,u)|^2=|x|^2+|u|^2$.

With this notation, the special problem can be turned into a special case of the General Problem by a change of variables:
\begin{align*}
	y_t&=\averageE\vect\left((x_t,x_{t-1},u_{t-1})x_t^\top\right)\\
	v_t&=\averageE\vect\left((x_t,u_t)u_t^\top\right)
\end{align*} 
The constraints $y_t\in\mathcal{X}$, $(y_t,v_t)\in\mathcal{U}$ and $(y_t,v_t,y_{t+1})\in\mathcal{V}$ correspond to positive semi-definiteness of the matrices
\begin{align*}
	&\averageE\left(x_tx_t^\top\right)\\
	&\averageE(x_t,u_t)(x_t,u_t)^\top\\
	&\averageE(x_t,u_t,x_{t+1})(x_t,u_t,x_{t+1})^\top.
\end{align*}

\section{Background}

The history of adaptive control dates back at least to the 1950s and the control of aircraft autopilots. The development was further stimulated by techniques for computer control and system identification in process industry. Following the landmark paper \cite{ast+wit73}, a surge of research activity during the 1970s derived conditions for convergence, stability, robustness and performance under various assumptions. Altogether, the subject has a rich history documented in numerous textbooks, such as \cite{goodwin2014adaptive}. 

The closely related term dual control was introduced by \cite{feldbaum1960dual} to describe the tradeoff between short term control objectives and actions to promote learning. 
This tradeoff is fundamentally important and has been studied extensively \cite{86helmersson+,88bernhardsson,wittenmark1995adaptive,filatov2000survey}. 
Dual control has also recently received renewed attention in connection machine learning. See \cite{mesbah2018stochastic}. 

In this paper, the focus is on worst-case models for disturbances and uncertain parameters, as discussed in \cite{cusumano1988nonlinear,sun1987theory,vinnicombe2004examples,2004megretskinonlinear} and more recently in \cite{rantzer2021minimax,cederberg2022synthesis,kjellqvist2022minimax,Rantzer2025acc}. 

Research on positive systems has a long history, dating back to the Luenberger's book \cite{Luenberger}, who introduced the concept together with a number of motivating examples. The use of nonnegative matrices when studying large scale systems was recognised by, for instance, \cite{Moylan+78} and \cite{sandell1978survey}. 
A more complete account on these subjects is given in the surveys \cite{BigSurvey,RantzerValcher2021}.

Control on more general positive cones has received relatively little attention so far, with the important exceptions of ensemble control \cite{li2010ensemble}, optimal transport theory \cite{chen2021optimal} and distributionally robust control \cite{van2015distributionally}. Our preliminary work \cite{pates2024optimal} can be viewed as a first step towards a unifying framework for control on positive cones.

\section{Solution to the General Problem.}

\begin{thm}
	Consider proper cones $\mathcal{X}\subseteq\realR^{N_y}$,  $\mathcal{U}\subseteq\realR^{N_y+N_v}$ and $\mathcal{V}, \mathcal{S}\subseteq\realR^{2N_y+N_v}$. Suppose that for every $y\in\mathcal{X}$ there exists $v$ with $(y,v)\in\mathcal{U}$. Assume also that the cones $\mathcal{U}$ and $\mathcal{V}$ have the property that
	\begin{align}
		\maxxi_{(y,v,\xi)\in\mathcal{V}}\minnu_{(\xi,\nu)\in\mathcal{U}}
		q^\top\big(z+(y,v,\xi),\xi,\nu\big)
	\label{eqn:linear}
	\end{align}
	is a linear function of $(z,y,v)$ whenever finite. Then the following conditions are equivalent:
	\begin{itemize}
		\item[\itemi]There exists a control law 
		$v_t=\mu_t(y_0,\ldots,y_t)$ 
		such that $(y_t,v_t)\in\mathcal{U}$ and \eqref{eqn:generalcost} holds
		for all sequences $y_0,y_1,\ldots$ with $y_0=0$ and $(y_t,v_t,y_{t+1})\in\mathcal{V}$ for all $t$.\smallskip
		
		\item[\ii]There exists a proper cone $\mathcal{Q}\subseteq\realR^{3N_y+2N_v}$ that contains $\mathcal{S}\times\{0,0\}$ and for all $y,v,z$ satisfies
		{\small\begin{align}
			\max_{q\in\mathcal{Q}}q^\top(z,y,v)
			&=\max_{q\in\mathcal{Q}}\maxxi_{(y,v,\xi)\in\mathcal{V}}\minnu_{(\xi,\nu)\in\mathcal{U}}
			q^\top\big(z+(y,v,\xi),\xi,\nu\big).
		\label{eqn:Bellman}
		\end{align}
		}
	\end{itemize}
	Moreover, if {\ii} holds, then {\itemi} holds with the feedback law
	{\small\begin{align}
		v_t&:=
		\arg\minv_{(y_t,v)\in\mathcal{U}}
		\max_{q\in\mathcal{Q}}q^\top\left(\sum_{\tau=0}^{t-1}(y_\tau,v_\tau,y_{\tau+1}),y_t,v\right).
	\label{eqn:controller}
	\end{align}}
\label{thm:general}
\end{thm}

\medskip

\begin{pf}
	To get that {\ii} $\Rightarrow$ {\itemi}, put
	\begin{align}
		z_t=\sum_{\tau=0}^{t-1}(y_\tau,v_\tau,y_{\tau+1})
	\label{eqn:z}
	\end{align} 
	and apply the controller \eqref{eqn:controller} to get
	\begin{align*}
		&\max_{q\in\mathcal{Q}}q^\top\left(\sum_{\tau=0}^{t-1}(y_\tau,v_\tau,y_{\tau+1}),y_t,v_t\right)\\
		&=^\eqref{eqn:z}\max_{q\in\mathcal{Q}}q^\top\left(z_t,y_t,v_t\right)\\
		&=^\eqref{eqn:Bellman}\max_{q\in\mathcal{Q}}\maxxi_{(y_t,v_t,\xi)\in\mathcal{V}}\minnu_{(\xi,\nu)\in\mathcal{U}}
q^\top\big(z_t+(y_t,v_t,\xi),\xi,\nu\big)\\
		&=\maxxi_{(y_t,v_t,\xi)\in\mathcal{V}}\max_{q\in\mathcal{Q}}
\minnu_{(\xi,\nu)\in\mathcal{U}}q^\top\big(z_t+(y_t,v_t,\xi),\xi,\nu\big)\\
		&=\maxxi_{(y_t,v_t,\xi)\in\mathcal{V}}\minnu_{(\xi,\nu)\in\mathcal{U}}
\max_{q\in\mathcal{Q}}q^\top\big(z_t+(y_t,v_t,\xi),\xi,\nu\big)\\
		&\ge\minnu_{(y_{t+1},\nu)\in\mathcal{U}}
\max_{q\in\mathcal{Q}}q^\top\big(z_t+(y_t,v_t,y_{t+1}),y_{t+1},\nu\big)\\
		&=^\eqref{eqn:controller}\max_{q\in\mathcal{Q}}q^\top\big(z_t+(y_t,v_t,y_{t+1}),y_{t+1},v_{t+1}\big)\\
		&=^\eqref{eqn:z}\max_{q\in\mathcal{Q}}q^\top\left(\sum_{\tau=0}^t(y_\tau,v_\tau,y_{\tau+1}),y_{t+1},v_{t+1}\right),
	\end{align*}
	where the fourth equality follows from von Neumann's minimax theorem.
	Recursive application of the resulting inequality starting from $y_0=0$ gives {\itemi}. 
	
	For the opposite implication, define $\mathcal{Q}_0:=\mathcal{S}\times\{0,0\}$ and generate an expanding sequence
	\begin{align*}
		\mathcal{Q}_0\subseteq
		\mathcal{Q}_1\subseteq
		\mathcal{Q}_2\subseteq\ldots
	\end{align*}
	recursively as follows. Assume that \eqref{eqn:Bellman} fails for $\mathcal{Q}=\mathcal{Q}_k$. This means that the value of \eqref{eqn:linear} is a linear function $q_k^\top(z,y,v)$ with larger value than $\max_{q\in\mathcal{Q}}q^\top(z,y,v)$ at some point. Then define
	\begin{align*}
		\mathcal{Q}_{k+1}&:=\conv(\mathcal{Q}_k\cup \{q_k\}).
	\end{align*}
	Assuming that {\itemi} holds, the recursion must converge to a limit set that satifies {\ii}.
\end{pf}

\goodbreak

\section{The Linear Quadratic Special Case}
\label{sec:solution}

We will now address the important instance of the 
Special Problem defined in section~\ref{sec:intro}, devoted to minimax optimal adaptive control for linear systems with quadratic cost. To do this, we need some notation. Given $\gamma>0$ and a set $\mathcal{M}$ of triples $(A,B,M)$, define $\mathcal{S}$ as the set of matrices 
\begin{align*}
  S&=\diag\{M,0\}-\gamma^2\begin{bmatrix}A&B&-I\end{bmatrix}^\top
  \begin{bmatrix}A&B&-I\end{bmatrix}
\end{align*}
generated by $(A,B,M)\in\mathcal{M}$.
In this setting, the solution to the general problem can be made more concrete and the following theorem can be stated:

\goodbreak

\begin{thm}
	Given $\mathcal{M}$ and $\mathcal{S}$, the following are equivalent:
	\begin{itemize}
		\item[\itemi]A control law 
		$u_t=\mu_t(x_0,\ldots,x_t)$ exists, such that the gain bound \eqref{eqn:gainbound} holds for all solutions to \eqref{eqn:state} with $(A,B,M)\in\mathcal{M}$.\smallskip
		
		\item[\ii]A proper cone $\mathcal{Q}\subseteq\realR^{(2n+m)\times(2n+m)}\times\realR^{(n+m)\times (n+m)}$ exists, that contains $\mathcal{S}\times\{0\}$ and satisfies
		{\small\begin{align}
				&\max_{(S,Q)\in\mathcal{Q}}\left(\langle S,Z\rangle +|(x,u)|^2_Q\right)\label{eqn:BellmanQ}\\
				&=\max_{\zeta}\min_\nu\max_{(S,Q)\in\mathcal{Q}}
				\averageE\left(\langle S,Z\rangle+|(x,u,\zeta)|^2_S+|(\zeta,\nu)|^2_Q\right)
				\notag
			\end{align}
		}for all $x,u,Z$.
	\end{itemize}
	Moreover, if {\ii} holds, then {\itemi} holds with the feedback law
	{\small\begin{align*}
			u_t&:=
			\arg\min_u
			\max_{(S,Q)\in\mathcal{Q}}\averageE
			\left(\langle S,Z_t\rangle+|(x_t,u_t)|^2_Q\right).
	\end{align*}}
	\label{thm:LQ}
\end{thm}

\begin{rmk}
	It should be noted that without randomization, the convex-concavity needed for application of van Neumann's minimax theorem would not be valid. This is what activates the input to get exploration when available data is insufficient. 
\end{rmk}

\begin{pf}
	The proof is analogous to Theorem~\ref{thm:general}.
\end{pf}

\section{Example with Explicit Solution}

For any $\alpha\in\realR_+$, define $\gamma_\alpha:=\alpha+\sqrt{1+\alpha^2}$ and let $\mathcal{M}_\alpha$ be the set of triples $(A,B,M)$ such that $A^\top A\preceq\alpha^2I$ and $B^\top B=I$, while 
\begin{align*}
	|x|^2&\ge\min_u\max_\zeta\left(|(x,u)|^2_M-\gamma_\alpha^2|Ax+Bu-\zeta|^2+|\zeta|^2\right).
\end{align*}

Then there exists a control law $\mu$ such that $$\averageE\sum_{t=0}^{T-1}|(x_t,u_t)|_M^2
\le \gamma_\alpha^2\sum_{t=0}^{T-1}|w_t|^2$$
whenever $x_{t+1}=Ax_t+Bu_t+w_t$, $u_t=\mu_t(x_0,\ldots,x_t)$, $x_0=0$ and $(A,B,M)\in\mathcal{M}_\alpha$. 

This follows from Theorem~\ref{thm:LQ} by explicit construction of $\mathcal{Q}$. For each triple $(A_i,B_i,M_i)\in\mathcal{M}$, define
{\small\begin{align*}
	S_{\pm i}&:=\diag\{M_i,0\}-\gamma^2\begin{bmatrix}A_i&\pm B_i&-I\end{bmatrix}^\top\begin{bmatrix}A_i&\pm B_i&-I\end{bmatrix}\\
	S_{0i}&:=(S_i+S_{-i})/2\\
	Q_i&:=M_i+\begin{bmatrix}A_i&B_i\end{bmatrix}^\top
	\begin{bmatrix}A_i&B_i\end{bmatrix}/(1-\gamma_\alpha^{-2})\\
	Q_0&:=M_i+\frac{\gamma_\alpha^2+1}{1-\gamma_\alpha^{-2}}\begin{bmatrix}A_i^\top A_i&0\\0&0\end{bmatrix}
	-\gamma_\alpha^2\begin{bmatrix}0&0\\0&B_i^\top B_i\end{bmatrix}\\
		\mathcal{Q}&:=\conv\cup_i\{(S_i,Q_i),(S_{0i},Q_{0i})\}.
	\end{align*}
}It is straightforward to verify that 
{\small\begin{align*}
	&\max_{j\in\{\pm 1,0\}}\max_i\left(\langle S_{ji},Z\rangle +|(x,u)|^2_{Q_{ji}}\right)\\
	&=\max_{\zeta}\min_\nu\max_{j\in\{\pm 1,0\}}\max_i
	\averageE\left(\langle S_{ji},Z\rangle+|(x,u,\zeta)|^2_{S_{ji}}+|(\zeta,\nu)|^2_{Q_{ji}}\right)
\end{align*}
}so equation \eqref{eqn:BellmanQ} follows.

\section{Acknowledgements}
The author is a member of the Excellence Center ELLIIT and Wallenberg AI, Autonomous Systems and Software Program (WASP). Support was received from the European Research Council (Advanced Grant 834142).


\end{document}